\documentclass[11pt, reqno]{amsart}
\usepackage{mathrsfs}
\usepackage[active]{srcltx}
\usepackage{mathrsfs}

\usepackage{color}
\usepackage[unicode]{hyperref}
\hypersetup{
	colorlinks = true,%
	citecolor = [rgb]{0.0,0.0,0.9},
	filecolor=black,%
	linkcolor = [rgb]{0.65,0.0,0.0},%
	anchorcolor = red,
	pagecolor = red,
	urlcolor= [rgb]{0.65,0.0,0.0}
}

\numberwithin{equation}{section}

\newtheorem{proposition}{Proposition}[section]
\newtheorem{theorem}{Theorem}[section]

\newtheorem{remark}{Remark}[section]

\newcommand{\ds}{\displaystyle}

\begin{document}
\title[Eigenvalues on negatively curved spaces]{New features of the first eigenvalue on negatively  curved spaces 
}

\vspace{-0.9cm}
\author{Alexandru Krist\'aly}
\address{Institute of Applied Mathematics, \'Obuda
	University, 
	Budapest, Hungary \& Department of Economics, Babe\c
	s-Bolyai University, Cluj-Napoca, Romania}
\email{kristaly.alexandru@nik.uni-obuda.hu; alexandrukristaly@yahoo.com 
}

\thanks{Research supported by the National Research, Development and Innovation Fund of Hungary, financed under the K$\_$18 funding scheme, Project No.  127926.}



\begin{abstract} 
	The paper is devoted to the study of fine properties of the first eigenvalue on negatively curved spaces. First, depending on the parity of the space dimension,  
	we provide asymptotically sharp harmonic-type expansions of the first eigenvalue for large geodesic balls in the model $n$-dimensional hyperbolic space, complementing the results of Borisov and Freitas (\textit{Comm. Anal. Geom.} 25:\ 507--544, 2017).
 We then give a  synthetic proof of Cheng's sharp eigenvalue comparison theorem in  metric measure spaces satisfying a 'negatively curved' Bishop-Gromov-type volume monotonicity hypothesis. 	
	As a byproduct,  we provide an example of simply connected, non-compact Finsler manifold with  constant negative flag curvature whose first eigenvalue is zero; this result is in a sharp contrast with its celebrated Riemannian counterpart due to McKean (\textit{J. Differential Geom.} 4:\ 359--366, 1970). Our proofs are based on specific properties of the Gaussian hypergeometric function combined with intrinsic aspects of the negatively curved smooth/non-smooth spaces.

\end{abstract}

\maketitle

\vspace{-0.0cm}
\section{Introduction and main results}
The goal of this paper is to establish new geometric properties encoded into the first eigenvalue on negatively curved (smooth or non-smooth) spaces.  In order to have a general geometric setting, we consider a (quasi)metric measure space  $(M,{d},\mu)$ with
a 
Borel measure $\mu,$ and let  ${\rm Lip}_0(\Omega)$ be the space of Lipschitz
functions with compact support on an open set $\Omega\subseteq M$. For  $u\in {\rm Lip}_0(\Omega)$,
let
\begin{equation}\label{local-constant}
|\nabla u|_{{{d}}}(x):=\limsup_{y\to x}
\frac{u(y)-u(x)}{{d}(x,y)}, \ x\in \Omega;
\end{equation}
note that $x\mapsto |\nabla u|_{d}(x)$ is Borel measurable on $\Omega$ and 
we may consider the \textit{fundamental frequency}  for $(\Omega,{d},\mu)$ defined by 
\begin{equation}\label{first-eigenvalue-general}
\lambda_{1,{d}}(\Omega):=\inf_{u\in {\rm Lip}_0(\Omega)\setminus \{0\}}\frac{\ds\int_{\Omega} |\nabla u|_{{{d}}}^2{\rm d}\mu}{\ds\int_{\Omega} u^2{\rm d}\mu}.
\end{equation}

\noindent In particular,  (\ref{first-eigenvalue-general}) corresponds to the \textit{first Dirichlet eigenvalue} of an open set $\Omega\subseteq M$ for the Laplace-Beltrami operator $-\Delta_g$  on a Riemannian manifold $(M,g)$ endowed with its usual canonical measure; a similar statement is also valid on Finsler manifolds with the Finsler-Laplace operator, see Ge and Shen \cite{GS}, Ohta and Sturm \cite{OS}.   

  On the one hand, when $(M,g)$ is a complete, simply  connected   $n$-dimensional   Riemannian  manifold  with sectional curvature bounded above by $-\kappa^2$ $(\kappa>0)$,   McKean \cite{McKean} proved in his celebrated paper that 
\begin{equation}\label{Mckeannn}
\lambda_{1,d_g}(M)\geq  \frac{(n-1)^2}{4}\kappa^2;
\end{equation} 
here, $d_g$ denotes the distance function on $(M,g)$.  
Moreover, in the  $n$-dimensional hyperbolic space $(\mathbb H^n_{-\kappa^2},g_h)$ of constant curvature $-\kappa^2$, the first eigenvalue has  the limiting property
\begin{equation}\label{elso-hatarertek}
\lim_{r\to \infty}\lambda_{1,d_h}(B_r^\kappa)=\lambda_{1,d_h}(\mathbb H^n_{-\kappa^2})=\frac{(n-1)^2}{4}\kappa^2,
\end{equation}
see Chavel \cite[p. 46]{Chavel} and Cheng and Yang \cite{Cheng-Yang}, where $B_r^\kappa$ and $d_h$ denote a geodesic ball of radius $r>0$ and the hyperbolic distance on $\mathbb H^n_{-\kappa^2}$, respectively.  

On the other hand, a consequence of the eigenvalue comparison theorem of Cheng \cite{Cheng} states that the hyperbolic space $\mathbb H^n_{-\kappa^2}$ has the greatest bottom of spectrum among all Riemannian manifolds with Ricci curvature bounded below by $-(n-1)\kappa^2$, i.e., 
\begin{equation}\label{Chenggggg}
\lambda_{1,d_g}(M)\leq  \frac{(n-1)^2}{4}\kappa^2.
\end{equation}
In the past half-century, McKean's and Cheng's results have become a continuing source  of inspiration concerning the first eigenvalue problem on curved spaces; without seeking completeness, we recall the works of Carroll and Ratzkin \cite{CR}, Chavel \cite{Chavel}, Freitas,  Mao and Salavessa \cite{Freitasandal}, Gage \cite{Gage}, Hurtado,  Markvorsen and Palmer \cite{HMP},    Li and Wang \cite{Li-Wang-1, Li-Wang-2}, Lott \cite{Lott}, Mao \cite{Mao}, Pinsky \cite{Pinsky, Pinsky-2} and  Yau \cite{Yau}, where various estimates and rigidity results concerning the equality in (\ref{Chenggggg}) are established.  

In view of (\ref{elso-hatarertek}) and (\ref{Chenggggg}), a considerable interest has been attracted to estimate the first eigenvalue of geodesic balls of  
 $(\mathbb H^n_{-\kappa^2},g_h)$ by means of elementary expressions. The most classical result states that  for every $n\geq 2$ one has
\begin{equation}\label{nullaban-aszimptotikus}
\lambda_{1,d_h}(B_r^\kappa)\sim\frac{j^2_{\frac{n}{2}-1,1}}{r^2}+\frac{n(n-1)}{6}\kappa^2\ \ {\rm as} \ \ r\to 0,
\end{equation}
see  Chavel \cite[p. 318]{Chavel}, where $j_{\frac{n}{2}-1,1}$ is the first positive zero of the Bessel function of first kind $J_{\frac{n}{2}-1}. $

In a recent result of Borisov and Freitas \cite[Theorem 3.3]{BF} the following two-sided estimate can be found for $\kappa=1$ (when we use the notation $B_r$ instead of $B_r^\kappa$) and arbitrary $r>0$:  
\begin{equation}\label{Borisov-2}
 \frac{j_{0,1}^2}{r^2}+\frac{1}{4}\left[\frac{1}{r^2}-\frac{1}{\sinh^2(r)}+1\right]\leq \lambda_{1,d_h}(B_r)\leq \frac{j_{0,1}^2}{r^2}+\frac{1}{3},\ \ \ \ \ \ \ \ \ \ \ \ \ \ \ \ \ \ \ \ \ \ \  \ \ \ \ \ \ \ \ \  \ \ n=2,
\end{equation}
$$ \frac{j_{\frac{n}{2}-1,1}^2}{r^2}+\frac{n(n-1)}{6}\leq \lambda_{1,d_h}(B_r)\leq \frac{j_{\frac{n}{2}-1,1}^2}{r^2}\ \ \ \ \ \ \ \ \ \ \ \ \ \ \ \ \ \ \ \ \ \ \ \ \ \ \ \  \ \ \ \ \ \ \ \ \ \ \ \ \ \ \ \ \ \ \ \ \ \ \ \ \ \ \ \ \ \ \ \ \ \ \ \ $$
\begin{eqnarray}\label{Borisov-n}
\ \ \ \ \ \ \ \ \ \ \ \  \ \ \ \ \ \  \ \ \ \ \ \ \ \ \ \ \ \ \ \ \ \ \ \ \ \ \ \ \ \ \ \  +\frac{(n-1)^2}{4}+\frac{(n-1)(n-3)}{4}\left[\frac{1}{\sinh^2(r)}-\frac{1}{r^2}\right],\ \ \ n\geq 3.
\end{eqnarray}
Since $j_{\frac{1}{2},1}=\pi,$ the estimates (\ref{Borisov-n}) spectacularly give in $3$-dimension the relation $\lambda_{1,d_h}(B_r)=1+\frac{\pi^2}{r^2}$ for every $r>0.$ 
We notice that the two-sided estimates (\ref{Borisov-2}) and (\ref{Borisov-n}) are asymptotically sharp for \textit{small} radii, i.e., the latter relations imply  (\ref{nullaban-aszimptotikus})  at once.  
However, apart from the case $n= 3$,  the estimates (\ref{Borisov-2}) and (\ref{Borisov-n}) are \textit{not}  asymptotically sharp whenever $r\to \infty$, see (\ref{elso-hatarertek}); only the lower bound in (\ref{Borisov-2}) and the upper bound in (\ref{Borisov-n}) behave properly, having their limit $\frac{(n-1)^2}{4}$ as $r\to \infty.$ 

Another estimate of $\lambda_{1,d_h}(B_r^\kappa)$ -- comparable to (\ref{Borisov-2}) and (\ref{Borisov-n}) -- which behaves  accurately for $r>0$ \textit{large} is provided by  Savo \cite[Theorem 5.6 (i)]{Savo} (see also Artamoshin \cite{Artamosin}), stating that for every $r>0$:
$$\frac{(n-1)^2}{4}\kappa^2+\frac{\pi^2}{r^2}-\frac{4\pi^2}{(n-1)r^3}\leq \lambda_{1,d_h}(B_r^\kappa)\leq  \frac{(n-1)^2}{4}\kappa^2+\frac{\pi^2}{r^2}+\frac{C}{r^3},$$
where $$C=\frac{\pi^2(n^2-1)}{2}\int_0^\infty \frac{s^2}{\sinh^2(s)}{\rm d}s.$$
In particular, one clearly has that 
\begin{equation}\label{SAvo}
\lambda_{1,d_h}(B_r^\kappa)= \frac{(n-1)^2}{4}\kappa^2+\frac{\pi^2}{r^2}+O\left({r^{-3}}\right)\ \ {\rm as}\ \ r\to \infty.
\end{equation}

Our first main result gives not only a more precise asymptotic behavior  than  (\ref{SAvo}) for \textit{large} radii $r>0$ (see also Cheng \cite[p. 294]{Cheng} and Borisov and Freitas \cite{BF})   but also provides a generic iterative method to compute/estimate $\lambda_{1,d_h}(B_r^\kappa)$ with respect to the space dimension (applicable mainly in the odd-dimensional case). 
In  order to state  our result, we introduce the  auxiliary functions
$$S_1(\gamma,x)=\frac{\sin(\gamma x)}{\gamma \sinh(x)}\ \ {\rm and}\ \ S_{k}(\gamma,x)=\frac{\frac{\partial S_{k-1}}{\partial x}(\gamma,x)}{\sinh(x)},\ k\geq 2,\ \gamma,x>0. $$

\begin{theorem}\label{fotetel-hiperbolikus} 
	Let $n\geq 2$ and $\kappa,r>0$. 
	%
	\begin{enumerate}
		\item[(i)] $(${\rm Odd-dimensional case}$)$	If $n=2l+1$ $(l\in \mathbb N),$ then $$\lambda_{1,d_h}(B_r^\kappa)=\frac{(n-1)^2}{4}\kappa^2+{\alpha^2},$$ where $\alpha=\alpha(\kappa,r,n)$ is the smallest positive solution to the transcendental equation  $S_l(\frac{\alpha}{\kappa},\kappa r)=0;$ in addition, for every $l\geq 2$, 
		$$\lambda_{1,d_h}(B_r^\kappa)\sim  \frac{(n-1)^2}{4}\kappa^2+\frac{\pi^2}{r^2}\left[1+\frac{1}{r}\left(1+\frac{1}{2}+...+
		\frac{1}{l-1}\right)\right]^2\ \ {as}\ \ r\to \infty.$$
		
		\item[(ii)] $(${\rm Even-dimensional case}$)$ If $n=2l$ $(l\in \mathbb N),$ then
		$$\lambda_{1,d_h}(B_r^\kappa)\sim  \frac{(n-1)^2}{4}\kappa^2+\frac{\pi^2}{r^2}\left[1+\frac{2}{r}\left(1+\frac{1}{3}+...+\frac{1}{2l-3}-\ln 2\right)\right]^2\ \ {as}\ \ r\to \infty.$$
		For $n=2$, the interior parenthesis reads as $-\ln 2.$
	\end{enumerate}
\end{theorem}


\begin{remark}\rm
	(i) In the particular case when $n=3$ (and $\kappa,r>0$ are fixed), the  
	smallest positive solution to the transcendental equation  $S_1(\frac{\alpha}{\kappa},\kappa r)=0$ is precisely $\alpha=\frac{\pi}{r};$ thus,  $\lambda_{1,d_h}(B_r^\kappa)=\kappa^2+\frac{\pi^2}{r^2}$ for every $r>0.$ This result  (for $n=3$)  coincides with the one of  Borisov and Freitas \cite{BF} and Savo \cite[Theorem 5.6 (ii)]{Savo},  where variational Hadamard-type formula and fine analysis on differential forms have been employed, respectively.  
	When $n\neq 3$, the above expressions provide the first four terms in the expansion of $\lambda_{1,d_h}(B_r^\kappa)$ for large $r>0$. Moreover, due to the alternating harmonic series $1-\frac{1}{2}+\frac{1}{3}-...=\ln 2,$ we have   $\left(1+\frac{1}{2}+...+
	\frac{1}{l-1}\right)\sim 2\left(1+\frac{1}{3}+...+\frac{1}{2l-3}-\ln 2\right)$ as $l\to \infty.$
	Thus, when the dimension is large enough (no matter on its parity), the lower-order terms in Theorem \ref{fotetel-hiperbolikus} (i) and (ii)  have similar asymptotic behavior.
	
	(ii) The transcendental equation $S_l(\frac{\alpha}{\kappa},\kappa r)=0$  in the odd-dimensional case $n=2l+1$ can be used to establish the asymptotically sharp form of $\lambda_{1,d_h}(B_r^\kappa)$ not only for large $r>0$, but also when $r\to 0$,   see (\ref{nullaban-aszimptotikus}); we exemplify this approach in dimension $n=5,$ see Remark \ref{remark-n-5} (i).

		(iii) The proof of Theorem \ref{fotetel-hiperbolikus} -- which is splitted according to the parity of the space dimension -- is based on a careful analysis of the Gaussian hypergeometric function whose first zero (with respect to certain parameter) is exactly the first eigenvalue $\lambda_{1,d_h}(B_r^\kappa)$,  see Section \ref{section-hyperbolic}. 
	
\end{remark}


Closely related to Theorem \ref{fotetel-hiperbolikus} -- where the Gaussian hypergeometric function appears as an extremal function on $B_r^\kappa\subset \mathbb H^n_{-\kappa^2}$ for $\lambda_{1,d_h}(B_r^\kappa)$ --  we establish a Cheng-type comparison result on metric measure spaces having a negatively curved character.  
 To be more precise, let $(M,d,\mu)$ be a (quasi)metric measure space with a strictly positive Borel measure $\mu$, and $x_0\in M$,  $\kappa>0$ and $n\in \mathbb N$ ($n\geq 2$) be fixed. We  assume first that small metric spheres in $M$ with center $x_0$ are comparable with their Euclidean counterparts; namely, we require the \textit{local density} assumption
\begin{description}
	\item[{\rm $({\bf D})^n_{x_0}$}]
	$\displaystyle\liminf_{\rho\to
		0}\frac{{\sf A}^\mu_{\rho}(x_0)}{n\omega_n\rho^{n-1}}=1$,
\end{description}
where 
\begin{equation}\label{area-keplet}
{\sf A}^\mu_{\rho}(x_0):=\frac{\rm d}{{\rm d}\rho}{\mu}(B_{\rho}(x_0))=\limsup_{\delta\to 0}\frac{\mu\left(B_{\rho+\delta}(x_0)\setminus B_{\rho}(x_0)\right)}{\delta}
\end{equation}
 denotes the induced $\mu$-area of the metric sphere $\partial B_{\rho}(x_0)$. Here,   $B_{\rho}(x_0)=\{y\in M:d(x_0,y)<\rho\},$ and $\omega_n$ is the volume of the $n$-dimensional Euclidean unit ball. Moreover, we introduce the following \textit{Bishop-Gromov-type volume monotonicity} hypothesis on the measure $\mu$: 
 \begin{description}
 	\item[{\rm $({\bf BG})^{n,\kappa}_{x_0}$}] the function 
 	$\rho\mapsto \frac{{\sf A}^\mu_{\rho}(x_0)}{\sinh^{n-1}(\kappa \rho)}$  is non-increasing on $(0,\infty)$. 
 \end{description}
  For further use, $V_\rho^\kappa$ stands for the hyperbolic volume of the  ball $B_\rho^\kappa\subset \mathbb H^n_{-\kappa^2}$.  
 

A sharp non-smooth  eigenvalue comparison principle of Cheng \cite{Cheng} (see also Hurtado, Markvorsen and Palmer \cite[Theorem E]{HMP})  reads as follows.

\begin{theorem}\label{fotetel-CD}  
	Let $(M,d,\mu)$ be a proper $($quasi$)$metric measure space, and assume the hypotheses {\rm $({\bf D})^n_{x_0}$} and $({\bf BG})^{n,\kappa}_{x_0}$ hold for some $x_0\in M$,  $\kappa>0$ and $n\in \mathbb N$ $(n\geq 2)$. 
	If  $r>0$ is fixed, then  
	\begin{equation}\label{2-becsles}
	\lambda_{1,d}(B_r(x_0))\leq \lambda_{1,d_h}(B_r^\kappa).
	\end{equation}
	Moreover,  if equality holds in  {\rm (\ref{2-becsles})} then $\mu(B_\rho(x_0))=V_\rho^\kappa$ for every $0<\rho<r$. 
\end{theorem}

\begin{remark}\rm

	(i) Hypothesis $({\bf BG})^{n,\kappa}_{x_0}$ is related to negative curvature; indeed,   $({\bf BG})^{n,\kappa}_{x}$  trivially holds on the hyperbolic space $\mathbb H^n_{-\kappa^2}$ for every $x\in \mathbb H^n_{-\kappa^2}$,  the function appearing in the hypothesis being constant. More generally, if a metric measure space $(M,d,\mu)$  satisfies the curvature-dimension condition
	${\sf CD}(-(n-1)\kappa^2,n)$ of Lott-Sturm-Villani 
	 for some $\kappa>0$ and $n\in \mathbb N$,  then the generalized Bishop-Gromov comparison principle states the validity of $({\bf BG})^{n,\kappa}_{x}$ for every $x\in M$, see Lott and Villani \cite{LV} and Sturm \cite{Sturm-Acta-II}. 
	However, there are metric measure spaces verifying $({\bf BG})^{n,\kappa}_{x_0}$ and failing   ${\sf CD}(-(n-1)\kappa^2,n)$ for every $\kappa>0$, see e.g.  the proof of  Theorem \ref{fotetel-Finsler} below. Another example is the Heisenberg group $(\mathbb H^m,d_{CC},\mathcal L^{2m+1})$ which verifies $({\bf BG})^{n,\kappa}_{x}$ for the homogeneous dimension $n=2m+2$  of  $ \mathbb H^m$ and every $\kappa>0$, $x\in \mathbb H^m$, and failing 	${\sf CD}(K,N)$ for \textit{any} choice of $K,N\in \mathbb R$, see Juillet \cite{Juillet}. 
	
	(ii)  Letting $r\to \infty$ in (\ref{2-becsles}), a similar inequality as (\ref{Chenggggg}) can be deduced on metric measure spaces satisfying $({\bf BG})^{n,\kappa}_{x_0}$ for some $x_0\in M$ (or satisfying the ${\sf CD}(-(n-1)\kappa^2,n)$ condition).
	
	(iii) Theorem \ref{fotetel-CD} can be applied to state various Cheng-type  comparison results on Riemannian/Finsler manifolds with (weighted) Ricci curvature bounded below. Indeed, under certain assumptions on the measure $\mu$ on an $n$-dimensional Finsler manifold $(M,F)$, the lower bound for the weighted Ricci curvature is equivalent to the condition ${\sf CD}(-(n-1)\kappa^2,n)$ for some $\kappa>0$, see Ohta \cite{Ohta-CV}. Moreover, when $\mu$ is the Busemann-Hausdoff measure, the local density assumption $({\bf D})^n_{x_0}$ holds for every $x_0\in M$, see Shen \cite{Shen-volume}, and Krist\'aly and Ohta \cite{Kristaly-Ohta}. The equality  in (\ref{2-becsles}) implies certain (radial) curvature rigidity and isometry between $B_r(x_0)$ and $B_r^\kappa$, see Cheng \cite[Theorem 1.1]{Cheng} and Zhao and Shen \cite[Theorem 1.2]{Zhao-Shen}; the details are left to the interested reader. The above-sketched consequences of Theorem \ref{fotetel-CD} complement in several aspects the results concerning the first eigenvalue problem on \textit{compact} Riemannian/Finsler manifolds developed by Ge and Shen \cite{GS}, Lott \cite{Lott}, Shen,  Yuan and Zhao \cite{Shen-YZ}, Wang and  Xia \cite{Wang-Xia}, and Wu and Xin \cite{Wu-Xin}.  
	
	(iv)	We notice that Cheng's original technique for proving (\ref{Chenggggg}) -- where smooth objects are explored as Jacobi vector fields and further properties of the exponential map on  Riemannian manifolds with Ricci curvature bounded below -- cannot be applied in the non-smooth framework of Theorem \ref{fotetel-CD}. However, it turns out that a contradiction argument combined with fine properties of the Gaussian hypergeometric function and  the Bishop-Gromov-type volume monotonicity hypothesis provide an elegant proof of Theorem \ref{fotetel-CD}, see Section \ref{section-Cheng}.

\end{remark}

An unexpected byproduct of Theorem \ref{fotetel-CD} is the following result in the Finsler setting,  which is in a sharp contrast with the Riemannian McKean's lower estimate (\ref{Mckeannn}). 

 \begin{theorem}\label{fotetel-Finsler} For every integer $n\geq 2$ there is a non-compact, forward complete, simply connected $n$-dimensional Finsler manifold $(M,F)$  with cons\-tant negative flag curvature  such that
\begin{equation}\label{egyenlo-nulla}
\lambda_{1,d_F}(M)=0,
\end{equation}
 where $d_F$ is the induced distance function on $(M,F)$. 
 \end{theorem}
 
 \begin{remark}\rm 
   One of the simplest Finsler structures fulfilling the thesis of  Theorem \ref{fotetel-Finsler} is provided by the $n$-dimensional  Euclidean open unit ball $B^n$  $(n\geq 2)$ endowed with the \textit{Funk metric} $F$, see Section \ref{section-Funk}.  We note that $(B^n,F)$  is a non-reversible Finsler manifold with constant flag curvature $-\frac{1}{4}$, having also   negative  weighted $N$-Ricci curvature for every $N\in [n,\infty].$ Beside the direct consequence  of Theorem \ref{fotetel-CD}, we present two further independent proofs for (\ref{egyenlo-nulla}).
 \end{remark}

In Section \ref{section-preliminaries} we recall those notations and results which are indispensable in our study, as basic properties of the hyperbolic spaces and Gaussian hypergeometric function,  a useful change-of-variable formula on metric measure spaces, and some elements from  Finsler geometry. In Sections \ref{section-hyperbolic}, \ref{section-Cheng} and \ref{section-Funk} we prove Theorems \ref{fotetel-hiperbolikus}, \ref{fotetel-CD} and \ref{fotetel-Finsler}, respectively.

\section{Preliminaries}\label{section-preliminaries}


\subsection{Hyperbolic spaces.} Let $\kappa>0$.  For the $n$-dimensional hyperbolic space we use the
Poincar\'e ball model $\mathbb H^n_{-\kappa^2}=\{x\in \mathbb R^n:|x|<1\}$
endowed with the Riemannian metric $$g_{
	h}(x)=(g_{ij}(x))_{i,j={1,...,n}}=p^2_\kappa(x)\delta_{ij},$$ where
$p_\kappa(x)=\frac{2}{\kappa(1-|x|^2)}.$  $(\mathbb
H^n_{-\kappa^2},g_{h})$ is a Cartan-Hadamard manifold with constant
sectional curvature $-\kappa^2$; the canonical volume form,  hyperbolic gradient and hyperbolic Laplacian operator are 
\begin{equation}\label{volume-form-hyper}
{\text d}v_{g_{h}}(x) = p^n_\kappa(x) {\text d}x,\ \ \nabla_{g_{h}}u=\frac{\nabla u}{p_\kappa^2}\ \ {\rm and}\ \,\Delta_{g_h}u=p_\kappa^{-n}{\rm div}(p_\kappa^{n-2}\nabla u),
\end{equation}
respectively, where $\nabla$ and div denote the Euclidean gradient and divergence operator in $\mathbb R^n.$ 
The
hyperbolic distance is denoted by $d_h$; the distance between the origin and $x\in \mathbb H^n_{-\kappa^2}$ is
given by $$d_{h}(0,x)=\frac{1}{\kappa}\ln\left(\frac{1+|x|}{1-|x|}\right).$$
The volume of the geodesic ball $B_r^\kappa=\{x\in \mathbb H^n_{-\kappa^2}:d_{h}(0,x)<r\}$ is $$\ds V_r^\kappa=n\omega_n\int_0^r \left(\frac{\sinh(\kappa\rho)}{\kappa}\right)^{n-1}{\rm d}\rho.$$
When $\kappa=1$, we simply use the notation $B_r$ and $\mathbb H^n$ instead of $B_r^\kappa$ and $\mathbb H^n_{-\kappa^2}$, respectively.

%
%

\subsection{Gaussian hypergeometric function} 
For $a,b,c\in \mathbb C$ ($c\neq 0,-1,-2,...$) we recall the  Gaussian hypergeometric function defined by  
$${\bf F}(a,b;c;z)=1+\sum_{k\geq 1}\frac{(a)_k(b)_k}{(c)_k}\frac{z^k}{k!}$$
on the disc $|z|<1$ and extended by analytic continuation elsewhere, where $(a)_k=\frac{\Gamma(a+k)}{\Gamma(a)}$ denotes the Pochhammer symbol. 
The corresponding differential equation to $z\mapsto {\bf F}(a,b;c;z)$ is 
\begin{equation}\label{hyper-ODE}
z(1-z)w''(z)+(c-(a+b+1)z)w'(z)-abw(z)=0.
\end{equation}
We also recall the differentiation formula
\begin{equation}\label{F-differential}
\frac{\rm d}{{\rm d}z}{\bf F}(a,b;c;z)=\frac{ab}{c}{\bf F}(a+1,b+1;c+1;z).
\end{equation}

Let $n\geq 2$ be an integer,  $C>0$ be fixed, and consider the second-order ordinary differential equation 
\begin{equation}\label{ODE-0}
\left(\frac{\rho^{n-1}}{(1-\rho^2)^{n-2}}
f'(\rho)\right)'+C\frac{\rho^{n-1}}{(1-\rho^2)^{n}}f(\rho)=0,\ \ \rho\in \left[0,1\right),
\end{equation}
subject to the boundary  condition $f(0)=1.$
The following result will be crucial in our  investigations. 
\begin{proposition}\label{proposition-oscillatory}
	The differential equation {\rm (\ref{ODE-0})} is oscillatory $($i.e.,  its  solutions 
	have an infinite number of  zeros$)$ if and only if $C>{(n-1)^2}$. 
\end{proposition}

{\it Proof.}  First, we transform {\rm (\ref{ODE-0})} into certain oscillation-preserving equivalent forms which will be useful in the proof. 
 Let $t=\frac{\rho^2}{1-\rho^2}$ and consider the function $w(t)=f(\rho);$ then  (\ref{ODE-0}) is transformed into 
\begin{equation}\label{ODE-modified}
\left(p(t)
w'(t)\right)'+q(t)w(t)=0,\ \ t>0,
\end{equation}
where $p(t)=4(t(t+1))^\frac{n}{2}$ and $q(t)=C(t(t+1))^\frac{n-2}{2}$.  Expanding   (\ref{ODE-modified}),  we equivalently obtain 
\begin{equation}\label{ODE-modified-2}
t(t+1)w''(t)+n\left(t+\frac{1}{2}\right)w'(t)+\frac{C}{4}w(t)=0,\ \ t>0.
\end{equation}
The trivial change of variables $t=-z$ in (\ref{ODE-modified-2}) leads to a differential equation of the form (\ref{hyper-ODE}). Therefore, the non-singular solution of (\ref{ODE-modified-2}) (since $w(0)=1$) is given by 
$$
w(t)={\bf F}\left(\frac{n-1+\sqrt{(n-1)^2-C}}{2},\frac{n-1-\sqrt{(n-1)^2-C}}{2};\frac{n}{2};-t\right),\ \ t>0,$$
thus 
\begin{equation}\label{hypergeom-vegso-megoldas}
f(\rho)={\bf F}\left(\frac{n-1+\sqrt{(n-1)^2-C}}{2},\frac{n-1-\sqrt{(n-1)^2-C}}{2};\frac{n}{2};\frac{\rho^2}{\rho^2-1}\right),\ \ \rho\in [0,1),
\end{equation}
represents the solution of (\ref{ODE-0}) with $f(0)=1$. We now distinguish the following two cases. 

\textit{Case 1}: $C>{(n-1)^2}$.
Since $\ds\int_{\alpha}^{\infty}\frac{1}{p(t)}{\rm d}t<\infty$ for every $\alpha>0$, we may apply Sugie, Kita and Yamaoka \cite[Theorem 3.1]{SKY} (see also Hille \cite{Hille}), i.e.,  if $$p(t)q(t)\left(\ds\int_{t}^{\infty}\frac{1}{p(\tau)}{\rm d}\tau\right)^2\geq \frac{1}{4}\ \ {\rm for}\ \ t\gg 1,$$ then (\ref{ODE-modified}) is oscillatory. The latter requirement  trivially holds since $C>{(n-1)^2}$; thus (\ref{ODE-0}) is also oscillatory.

\textit{Case 2}: $C\leq {(n-1)^2}$. By (\ref{hypergeom-vegso-megoldas}) and the connection formula  (15.10.11) of \cite{Digital}, one has for every $\rho\in [0,1)$ that
$$f(\rho)=(1-\rho^2)^\frac{n-1+\sqrt{(n-1)^2-C}}{2}{\bf F}\left(\frac{n-1+\sqrt{(n-1)^2-C}}{2},\frac{1+\sqrt{(n-1)^2-C}}{2};\frac{n}{2};{\rho^2}\right)>0,$$ thus {\rm (\ref{ODE-0})} is non-oscillatory. 
 \hfill $\square$

%

\subsection{Change-of-variables formula.}

Let $({M},d,\mu)$ be a (quasi)metric measure space, i.e.,
$({M},d)$ is a complete separable (quasi)metric space and
$\mu$ is a locally finite measure on $M$ endowed with its
Borel $\sigma$-algebra. We assume that the measure
$\mu$ on $M$ is strictly positive, i.e.,
supp$[\mu]=M.$ 
Let $B_{\rho}(x)=\{y\in M:{{d}}(x,y)<\rho\}$. 

A useful change-of-variables formula on $(M,d,\mu)$ reads as follows. 



\begin{proposition}\label{layer-cake}
	Let $r>0$ and  $f:(0,r]\to \mathbb R$ be a non-increasing function such that $f(r)=0$,  $({M},d,\mu)$ be a $($quasi$)$metric measure space, and assume the hypothesis $({\bf BG})^{n,\kappa}_{x_0}$ holds for some $x_0\in M$,  $\kappa>0$ and $n\in \mathbb N$ $(n\geq 2)$.   Then  $$\int_{B_{r}(x_0)}f({d}(x_0,x)){\rm d}\mu(x)=\int_0^r {\sf A}^\mu_{\rho}(x_0)f(\rho){\rm d}\rho.$$
\end{proposition}

{\it Proof.} By hypothesis $({\bf BG})^{n,\kappa}_{x_0}$ and   Gromov's monotonicity result, see e.g.   Cheeger,  Gromov and Taylor \cite[p. 42]{CGT}, it follows that $\rho\mapsto \frac{\mu(B_{\rho}(x_0))}{V_\rho^\kappa}$ is non-increasing on $(0,\infty)$; in particular, $\rho\mapsto \mu(B_{\rho}(x_0))$ is differentiable a.e.\ on $[0,\infty)$. Let $l_0=\lim_{\rho\to 0} f(\rho).$ By the layer cake representation and the facts that  $f:(0,r]\to \mathbb R$ is non-increasing and $f(r)=0$, an integration by parts provides 
\begin{eqnarray*}
	\int_{B_{r}(x_0)}f({d}(x_0,x)){\rm d}\mu(x)&=&\int_0^{l_0}\mu(\left\{x\in B_{r}(x_0):f({d}(x_0,x))>t\right\}){\rm d}t\\&=&
	\int_r^0\mu( B_{\rho}(x_0))f'(\rho){\rm d}\rho \ \ \ \ \ \ \ \ \ \ \ \ \ \ \ [{\rm change\ of\ variables}\ t=f(\rho)]\\&=&\int_0^r \frac{\rm d}{{\rm d}\rho}\mu(B_{\rho}(x_0))f(\rho){\rm d}\rho,
\end{eqnarray*}
as we intended to prove. 
\hfill $\square$

\subsection{Finsler manifolds.} Let $M$ be a connected $n$-dimensional smooth manifold and $TM=\bigcup_{x \in M}T_{x} M $ be its tangent bundle. 
The pair $(M,F)$ is called a {\it Finsler manifold} if the
continuous function $F:TM\to [0,\infty)$ satisfies the
conditions:

(a) $F\in C^{\infty}(TM\setminus\{ 0 \});$

(b) $F(x,ty)=tF(x,y)$ for all $t\geq 0$ and $(x,y)\in TM;$

(c) $g_v=[g_{ij}(v)]:=\left[\frac12F^{2}%
(x,y)\right]_{y^{i}y^{j}}$ is positive definite for all $v=(x,y)\in
TM\setminus\{ 0 \},$\\
\noindent see Bao, Chern and Shen \cite{BCS}. 
If $F(x,ty)=|t|F(x,y)$ for every $t\in \mathbb R$ and $(x,y)\in TM,$ then $(M,F)$ is called reversible. If $g_{ij}(x)=g_{ij}(x,y)$ is
independent of $y$ then $(M,F)=(M,g)$ is a {Riemannian
	manifold}. 

For every
$(x,\alpha)\in T^*M$, the \textit{co-metric} (or, polar transform)
of $F$ is defined  by
\begin{equation}  \label{polar-transform}
F^*(x,\alpha)=\sup_{v\in T_xM\setminus
	\{0\}}\frac{\alpha(v)}{F(x,v)}.
\end{equation}

Unlike the Levi-Civita connection \index{Levi-Civita connection} in
the Riemannian case, there is no unique natural connection in the
Finsler geometry. Among these connections on the pull-back bundle
$\pi ^{*}TM,$ we choose a torsion free and almost metric-compatible
linear connection on $\pi ^{*}TM$, the so-called \textit{Chern
	connection}, see Bao, Chern and Shen \cite[Theorem 2.4.1]{BCS}.
The Chern connection induces on $\pi ^{*}TM$ the \textit{%
	curvature tensor} $R$.  The Finsler manifold is
{\it forward {\rm (resp.} bacward$)$ complete} if every
geodesic segment $\sigma:[0,a]\to M$ can be extended to $[0,\infty)$ (resp. $(-\infty,0]$).

Let $u,v\in T_xM$ be two non-collinear vectors and $\mathcal{S}=\mathrm{span}%
\{u,v\}\subset T_xM$. By means of the curvature tensor $R$, the
\textit{flag curvature} associated with the flag $\{\mathcal{S},v\}$ is 
\begin{equation}  \label{ref-flag}
\mathbf{K}(\mathcal{S};v) =%
\frac{g_v(R(U,V)V, U)}{g_v(V,V) g_v(U,U) - g_v^{2}(U,V)},
\end{equation}
where $U=(v;u),V=(v;v)\in \pi^*TM.$ If $(M,F)$ is Riemannian, the
flag
curvature reduces to the sectional curvature which depends only on $\mathcal{S}$.  

Take $v\in T_xM$ with $F(x,v)=1$ and
let $\{e_i\}_{i=1}^n$ with $e_n=v$ be an orthonormal basis of
$(T_xM,g_v)$ for $g_v$. Let $\mathcal S_i={\rm
	span}\{e_i,v\}$ for $i=1,...,n-1$. Then the {\it Ricci curvature} of
$v$ is defined by ${\rm Ric}(v):=\sum_{i=1}^{n-1}\textbf{K}(\mathcal S_i;v)$.

Let $\mu$ be a positive smooth measure on $(M,F)$. Given $v \in
T_xM \setminus \{0\}$, let $\sigma:(-\varepsilon,\varepsilon)
\to M$ be the geodesic with $\dot{\sigma}(0)=v$ and
decompose $\mu$ along $\sigma$ as
$\mu=e^{-\psi}\mathrm{vol}_{\dot{\sigma}}$, where
$\mathrm{vol}_{\dot{\sigma}}$ denotes the volume form of the
Riemannian structure $g_{\dot{\sigma}}$. For $N \in
[n,\infty]$, the \emph{weighted $N$-Ricci curvature} $\mathrm{Ric}_N$ is
defined by
\[ \mathrm{Ric}_N(v):=\mathrm{Ric}(v) +(\psi \circ \sigma)''(0)
-\frac{(\psi \circ \sigma)'(0)^2}{N-n}, \]
where the third term is understood as $0$ if $N=\infty$ or if $N=n$
with $(\psi \circ \sigma)'(0)=0$, and as $-\infty$ if $N=n$ with
$(\psi \circ \sigma)'(0) \neq 0$. 

Let $\sigma: [0,r]\to M$ be a piecewise smooth curve. The value $%
L_F(\sigma)= \displaystyle\int_{0}^{r} F(\sigma(t), \dot\sigma(t))\,{\text d}%
t $ denotes the \textit{integral length} of $\sigma.$ For
$x_1,x_2\in M$,
denote by $\Lambda(x_1,x_2)$ the set of all piecewise $C^{\infty}$ curves $%
\sigma:[0,r]\to M$ such that $\sigma(0)=x_1$ and $\sigma(r)=x_2$.
Define the \textit{distance function} $d_{F}: M\times M
\to[0,\infty)$ by
\begin{equation}  \label{quasi-metric}
d_{F}(x_1,x_2) = \inf_{\sigma\in\Lambda(x_1,x_2)} L_F(\sigma).
\end{equation}
If $u\in C^1(M)$,  on account of (\ref{local-constant}) we have 
\begin{equation}\label{osszefugges-F-gradiens}
|\nabla u|_{d_F}(x)=F^*(x,Du(x)), \ x\in M.
\end{equation}
In particular, if $x_0\in M$ is
fixed, then we have the \textit{eikonal equation}
\begin{equation}\label{tavolsag-derivalt}
F^*(x,D d_F(x_0,x))=1\ {\rm for\ a.e.}\ x\in
M.
\end{equation}

Let $\{{\partial}/{\partial x^i} \}_{i=1,...,n}$ be a local basis
for the tangent bundle $TM,$ and $\{\mathrm{d}x^i \}_{i=1,...,n}$ be
its dual basis
for $T^*M.$ Consider $\tilde B_x(1)=\{y=(y^i):F(x,y^i \partial/\partial x^i)< 1\}\subset \mathbb R^n$. The \textit{Busemann-Hausdorff volume form}  is defined by
\begin{equation}  \label{volume-form}
{\text d}v_F(x)=\sigma_F(x){\text
	d}x^1\wedge...\wedge {\text d}x^n,
\end{equation}
where $\sigma_F(x)=\frac{\omega_n}{|\tilde B_x(1)|}$. 


The {\it Legendre transform}
$J^*:T^*M\to TM$ associates to each element $\xi\in T_x^*M$ the
unique maximizer on $T_xM$ of the map $y\mapsto
\xi(y)-\frac{1}{2}F(x,y)^2$. The {\it gradient} of $u$ is defined by
${\nabla}_F u(x)=J^*(x,Du(x)).$
The {\it Finsler-Laplace operator} is given by 
$${\Delta}_F u={\rm div}_F({\nabla}_F u),$$ where 
div$_F(X)=\frac{1}{\sigma_F}\frac{\partial}{\partial
	x^i}(\sigma_F X^i)$ for some vector field $X$  on $M$, and $\sigma_F$ comes from (\ref{volume-form}).


\section{Proof of Theorem \ref{fotetel-hiperbolikus}}\label{section-hyperbolic}

By a density reason we have
for every $r>0$ that  $$\lambda_{1,{d_h}}(B_r^\kappa)=\inf_{u\in H_0^1(B_r^\kappa)\setminus \{0\}}\frac{\ds\int_{B_r^\kappa} |\nabla_{g_h} u|_{g_h}^2{\rm d}v_{g_h}}{\ds\int_{B_r^\kappa} u^2{\rm d}v_{g_h}},$$
where $H_0^1(B_r^\kappa)$ is the usual Sobolev space over the Riemannian manifold $(B_r^\kappa,g_h)$, see Hebey \cite{Hebey}. 
Standard arguments from calculus of variations -- based on the compactness of the embedding $H_{0}^1(B_r^\kappa) \hookrightarrow L^2(B_r^\kappa)$  and the convexity (thus, the sequentially weakly lower semicontinuity) of  $u\mapsto \ds\int_{B_r^\kappa} |\nabla_{g_h} u|_{g_h}^2{\rm d}v_{g_h}$ on $H_{0}^1(B_r^\kappa)$, -- imply the existence of a minimizer  for $\lambda_{1,d_h}(B_r^\kappa)$. The minimizer for  $\lambda_{1,d_h}(B_r^\kappa)$ can be assumed to be positive; moreover, a  convexity reason  shows that it is unique up to constant multiplication. Let $w^*:B_r^\kappa\to [0,\infty)$  be the positive minimizer for  $\lambda_{1,d_h}(B_r^\kappa)$. 
Moreover, one can deduce by a P\'olya-Szeg\H o-type inequality on  $(\mathbb H^n_{-\kappa^2},g_{h})$ (see 
Baernstein \cite{Baernstein}) that $w^*$ is radially symmetric.  In particular, by standard regularity and Euler-Lagrange equation, if $w^*(x)=f(|x|)$ with $f:[0,\tanh(\frac{\kappa r}{2}))\to [0,\infty)$  smooth enough,  we obtain 
\begin{equation}\label{ODE}
\left(\frac{\rho^{n-1}}{(1-\rho^2)^{n-2}}f'(\rho)\right)'+\frac{4\lambda_{1,d_h}(B_r^\kappa)}{\kappa^2}\frac{\rho^{n-1}}{(1-\rho^2)^{n}}f(\rho)=0,\ \ \rho\in \left[0,\tanh\left(\frac{\kappa r}{2}\right)\right),
\end{equation}
subject to the boundary condition $f(\tanh(\frac{\kappa r}{2}))=0$; the latter relation comes from the fact that $w^*$ vanishes on $\partial B_\kappa^r$. Since $\tanh(\frac{\kappa r}{2})<1$, in order to fulfill the boundary  condition, we need to guarantee the oscillatory behaviour of (\ref{ODE}). Due to Proposition \ref{proposition-oscillatory}, the latter statement is equivalent to  $\frac{4\lambda_{1,d_h}(B_r^\kappa)}{\kappa^2}>(n-1)^2$, which perfectly agrees with McKean's estimate (\ref{Mckeannn}).      Let
\begin{equation}\label{alpha-ertek}
\alpha^\kappa_{n,r}:=\sqrt{\frac{\lambda_{1,d_h}(B_r^\kappa)}{\kappa^2}-\frac{(n-1)^2}{4}}>0.
\end{equation}
By (\ref{hypergeom-vegso-megoldas}), it turns out that the solution of (\ref{ODE}) 
can be written into the form 
\begin{equation}\label{meg-ez-is-kell}
f(\rho)={\bf F}\left(\frac{n-1}{2}+ i\alpha^\kappa_{n,r},\frac{n-1}{2}-i\alpha^\kappa_{n,r};\frac{n}{2};\frac{\rho^2}{\rho^2-1}\right),\ \ \rho\in \left[0,\tanh\left(\frac{\kappa r}{2}\right)\right).
\end{equation}
%
The boundary  condition $f(\tanh(\frac{\kappa r}{2}))=0$ implies that  $\rho=\tanh(\frac{\kappa r}{2})$ is the first positive zero of (\ref{ODE}); therefore, the value of $\lambda_{1,d_h}(B_r^\kappa)$ is obtained as the smallest positive solution of 
\begin{equation}\label{utolso-az-elso-tetelnel}
{\bf F}\left(\frac{n-1}{2}+ i\alpha^\kappa_{n,r},\frac{n-1}{2}-i\alpha^\kappa_{n,r};\frac{n}{2};-\sinh^2\left(\frac{\kappa r}{2}\right)\right)=0,
\end{equation}
where $\alpha^\kappa_{n,r}$ is given in (\ref{alpha-ertek}). Having the (theoretical) value of $\lambda_{1,d_h}(B_r^\kappa)$, it turns out that 
\begin{equation}\label{w-csillag}
w^*(x)=f(|x|)=f\left(\tanh\left(\frac{\kappa d_{h}(0,x)}{2}\right)\right),\ \ x\in B_r^\kappa.
\end{equation}
By construction, $f(0)=1$ and  $f(\rho)>0$ for every $\rho\in \left[0,\tanh(\frac{\kappa r}{2})\right)$; moreover, a simple monotonicity reasoning based on (\ref{ODE}) shows that  $\rho\mapsto f(\rho)$ is decreasing on $\left[0,\tanh(\frac{\kappa r}{2})\right)$. 

\begin{remark}\rm 
	With the above notations, the form of the hyperbolic Laplacian in (\ref{volume-form-hyper}) shows that equation (\ref{ODE}) corresponds precisely to the eigenvalue problem 
	\begin{equation*}
\left \{
	\begin{array}{lll}
	 -\Delta_{g_h}w^*=\lambda_{1,d_h}(B_r^\kappa) w^* & {\rm in} &  B_r^\kappa, \\
	 w^*=0 & {\rm on} &  \partial B_r^\kappa.%
	\end{array}
	\right.
	\end{equation*}
\end{remark}

%

In the rest of this section we consider $\kappa=1$, the general case easily following by a scaling argument; we will use the notation $\alpha_{n,r}$ instead of $\alpha_{n,r}^\kappa.$\\ 

We now prove Theorem \ref{fotetel-hiperbolikus} by splitting our argument according to the parity of dimension.

\subsection{Odd-dimensional case.} This part is also divided into two sub-cases. 
\subsubsection{The case $n=3$}

First of all, we claim that  for every $\gamma>0$ and $\rho\in (0,1),$ one has the identity
\begin{equation}\label{azonossag}
{\bf F}\left(1+i\gamma,1-i\gamma;\frac{3}{2};\frac{\rho^2}{\rho^2-1}\right)=\frac{1-\rho^2}{2\gamma \rho}{\sin\left(2\gamma { \tanh}^{-1}(\rho)\right)}.
\end{equation}
To verify (\ref{azonossag}), we look for the solution of (\ref{ODE-0}) in the  form $f(\rho)=c_0\frac{1-\rho^2}{\rho}s(\rho)$ for some $c_0>0$ whenever $n=3$ and $C=4(\gamma^2+1)>4$.  Thus, a simple computation transforms (\ref{ODE-0}) into the equation $$\rho s''(\rho)-\frac{2\rho^2}{1-\rho^2}s'(\rho)+4\gamma^2\frac{\rho}{(1-\rho^2)^2}s(\rho)=0,\ \ \rho\in (0,1),$$
with the boundary  condition $s(0)=0.$ Now, if $\rho=\tanh(t)$ and $s(\rho)=w(t)$, the latter equation is transformed into 
$$w''(t)+4\gamma^2w(t)=0,\ \ t>0,$$
with the boundary  condition $w(0)=0;$ thus $w(t)=\sin(2\gamma t)$. Now, relation (\ref{azonossag}) follows by (\ref{hypergeom-vegso-megoldas}) and the fact that $f(0)=1$, thus $c_0=\frac{1}{2\gamma}$. 

Let us choose $\gamma:={\alpha_{3,r}}=\sqrt{{\lambda_{1,d_h}(B_r)}-1}$ and $\rho:=\tanh(\frac{ r}{2})$  in (\ref{azonossag}). Due to (\ref{utolso-az-elso-tetelnel}), the value $\lambda_{1,d_h}(B_r)$ is precisely the smallest positive solution of 
$\sin\left( r\sqrt{{\lambda_{1,d_h}(B_r)}-1}\right)=0,$ i.e.,  
$$\lambda_{1,d_h}(B_r)=1+\frac{\pi^2}{r^2}.$$



\subsubsection{The case $n=2l+1\geq 5$} The identity (\ref{azonossag}) can be equivalently written into the form
\begin{equation}\label{azonossag-2}
{\bf F}\left(1+i\gamma,1-i\gamma;\frac{3}{2};-\sinh^2(\frac{x}{2})\right)=\frac{\sin(\gamma x)}{\gamma \sinh(x)},\ \ \gamma,x>0. 
\end{equation}
For every $\gamma,x>0$, let us introduce the functions 
$$S_1(\gamma,x)=\frac{\sin(\gamma x)}{\gamma \sinh(x)}\ \ {\rm and}\ \ S_{k}(\gamma,x)=\frac{\frac{\partial S_{k-1}}{\partial x}(\gamma,x)}{\sinh(x)},\ k\geq 2. $$
By applying inductively the differentiation formula (\ref{F-differential}), we obtain for every $\gamma,x>0$ and integer $k\geq 1$ that 
\begin{equation}\label{iterative}
{\bf F}\left(k+i\gamma,k-i\gamma;\frac{2k+1}{2};-\sinh^2(\frac{x}{2})\right)=(-1)^{k-1}\frac{(2k-1)!!}{\prod_{j=1}^{k-1}(j^2+\gamma^2)}S_k(\gamma,x),
\end{equation}
where by convention the denominator at the right hand side is $1$ for $k=1.$
According to (\ref{alpha-ertek}), (\ref{utolso-az-elso-tetelnel}) and (\ref{iterative}), for $n=2l+1$, we have 
\begin{equation}\label{elso-sajat-alpha}
\lambda_{1,d_h}(B_r)=\frac{(n-1)^2}{4}+\alpha^2,
\end{equation}  
where $\alpha=\alpha(r,l)$ is the smallest positive solution to the transcendental equation $S_l(\alpha,r)=0$. Although no explicit solution $\alpha$ can be provided  to the latter equation, one can prove first that 
\begin{equation}\label{alpha-konv}
\alpha\sim\frac{\pi}{r}\ \ {\rm as}\ \ r\to \infty.
\end{equation}
Before to do this, let us observe that by the estimate (\ref{Borisov-n}) and (\ref{elso-sajat-alpha}), one has  $0<\alpha r\leq j_{l-\frac{1}{2},1}$ as $r\to \infty;$ thus, we may assume that $\alpha r\to \Phi$ as $r\to \infty$ for some $\Phi\in (0,j_{l-\frac{1}{2},1}].$ In particular, $\alpha\to 0$ as $r\to \infty$. We are going to prove that $\Phi=\pi,$ which completes (\ref{alpha-konv}).

As a model situation, let us consider some lower-dimensional cases. When $n=5$ (thus $l=2$), the equation $S_2(\alpha,r)=0$ is equivalent to 
\begin{equation}\label{S2-egyenlet}
\alpha\cos(\alpha r) \tanh(r)-\sin(\alpha r)=0;
\end{equation}
 taking the limit  $r\to \infty$, it follows that $\sin(\Phi)=0$, i.e., $\Phi=\pi,$ due to the minimality property of $\Phi>0$.  When $n=7$ (thus $l=3$), the equation $S_3(\alpha,r)=0$ is equivalent to $3\alpha \cos(\alpha r)\tanh(r)+\sin(\alpha r)[(\alpha^2+1)\tanh^2(r)-3]=0.$ A similar limiting reasoning as above gives $\sin(\Phi)=0$, thus $\Phi=\pi$.

In higher-dimensional cases the equation $S_l(\alpha,r)=0$ becomes more and more involved. In order to handle this generic case, let us observe that $\sinh(r)\sim e^r/2$ and $\cosh(r)\sim e^r/2$ as $r\to \infty,$ and for every smooth function $\Psi:(0,\infty)\to \mathbb R$, one has a stability property of differentiation  with respect to approximation of hyperbolic functions,  i.e.,  
$$\left.\frac{\rm d}{{\rm d}x}\frac{\Psi(x)}{\sinh(x)}\right|_{x=r}\sim 2(\Psi'(r)-\Psi(r))e^{-r}=2\left.\frac{\rm d}{{\rm d}x}({\Psi(x)}{e^{-x}})\right|_{x=r}\ \ {\rm as}\ \ r\to \infty.$$ 
Accordingly, in order to establish the asymptotic behavior of $\alpha$ with respect to $r$ when $r\to \infty,$ we may consider instead of  $S_l(\alpha,r)=0$  the approximation equation $\tilde S_l(\alpha,r)=0$, where 
$$\tilde S_1(\gamma,x)={\sin(\gamma x)}e^{-x}\ \ {\rm and}\ \ \tilde S_{k}(\gamma,x)={\frac{\partial\tilde S_{k-1}}{\partial x}(\gamma,x)}e^{-x},\ \ k\geq 2,\ \gamma,x>0. $$
By induction, one can easily prove that $$\tilde S_k(\gamma,x)=(P_k(\gamma)\cos(\gamma x)+Q_k(\gamma)\sin(\gamma x))e^{-kx},\ \ k\geq 1,\ \gamma,x>0,$$
where 
\begin{eqnarray}\label{P-Q}
\left\{
\begin{array}{lll}
P_{k+1}(\gamma)&=&\gamma Q_k(\gamma)-kP_k(\gamma), \\
Q_{k+1}(\gamma)&=&-\gamma P_k(\gamma)-kQ_k(\gamma),
\end{array}\right. k\geq 1,\ \gamma\geq 0,
\end{eqnarray}
%
%
%
%
and $P_1\equiv 0,$ $Q_1\equiv 1.$
We observe that $P_k(0)=0\neq Q_k(0)$ for every $k\geq 1.$ 
Now, by limiting in $\tilde S_l(\alpha,r)=0$ as $r\to \infty$ and taking into account that $\alpha\to 0$,  it turns out that $\sin(\Phi)=0$, i.e., $\Phi=\pi,$ which concludes the proof of (\ref{alpha-konv}). 

Let $n=2l+1$ with $l\geq 2$.  We prove that 
\begin{equation}\label{alpha-ujabb}
\alpha\sim\frac{\pi}{r}+\frac{c_l}{r^2}\ \ {\rm as}\ \ r\to \infty,
\end{equation}
for some $c_l\in \mathbb R$ which will be determined in the sequel. Plugging the latter form of $\alpha$ into the approximation equation $\tilde S_l(\alpha,r)=0$, one has approximately that
$$P_l\left(\frac{\pi}{r}+\frac{c_l}{r^2}\right)\cos\left(\frac{c_l}{r}\right)+Q_l\left(\frac{\pi}{r}+\frac{c_l}{r^2}\right)\sin\left(\frac{c_l}{r}\right)=0.$$
Multiplying the latter relation by $r>0$, letting $z=1/r$ and taking the limit when $r\to \infty$, it follows that
$c_l=-\lim_{z\to 0}\frac{P_l(\pi z)}{zQ_l(\pi z)}.$ Since 
$P_l(0)=0\neq Q_l(0)$, the Taylor expansion of $P_l$ and $Q_l$  gives that  
$c_l=-\frac{P_l'(0)}{Q_l(0)}\pi.$
By the second relation of (\ref{P-Q}), we directly obtain $Q_l(0)=(-1)^{l-1}(l-1)!$, while from the first relation we deduce the recurrence $P_{k+1}'(0)=Q_k(0)-kP_k'(0),$  $k\geq 1$. A simple reasoning implies that $P_l'(0)=(-1)^l(l-1)!(1+\frac{1}{2}+...+\frac{1}{l-1})$. Consequently, (\ref{alpha-ujabb}) follows since  $$c_l=\pi\left(1+\frac{1}{2}+...+\frac{1}{l-1}\right),\ \ l\geq2.$$


\subsection{Even-dimensional case.} Up to some technical differences, the structure of the proof is the same as in the odd-dimensional case. First of all,  one has  for every $\gamma,x>0$ that
\begin{equation}\label{spherical}
{\bf F}\left(\frac{1}{2}+i\gamma,\frac{1}{2}-i\gamma;1;-\sinh^2(\frac{x}{2})\right)={\bf P}_{-\frac{1}{2}+i\gamma}(\cosh(x)), 
\end{equation}
where ${\bf P}_{-\frac{1}{2}+i\gamma}$ denotes the spherical Legendre function, see Robin \cite{Robin},  Zhurina and Karmazina \cite{ZK}.  For every $\gamma,x>0$, we consider the functions 
$$\overline S_1(\gamma,x)={\bf P}_{-\frac{1}{2}+i\gamma}(\cosh(x))\ \ {\rm and}\ \ \overline S_{k}(\gamma,x)=\frac{\frac{\partial \overline S_{k-1}}{\partial x}(\gamma,x)}{\sinh(x)},\ k\geq 2. $$
By the differentiation formula (\ref{F-differential}) and (\ref{spherical}) we have for every $\gamma,x>0$ and integer $k\geq 1$ that 
\begin{equation}\label{iterative-2}
{\bf F}\left(\frac{2k-1}{2}+i\gamma,\frac{2k-1}{2}-i\gamma;k;-\sinh^2(\frac{x}{2})\right) =(-1)^{k-1}
\frac{2^{k-1}(k-1)!}{\prod_{j=1}^{k-1}(\frac{(2j-1)^2}{4}+\gamma^2)}\overline S_k(\gamma,x),
\end{equation}
where by convention the denominator at the right hand side is $1$ for $k=1.$ 

Let  $n=2l,$ $l\in \mathbb N$. Due to  (\ref{utolso-az-elso-tetelnel}) and (\ref{iterative-2}), we have 
\begin{equation}\label{elso-sajat-alpha-2}
\lambda_{1,d_h}(B_r)=\frac{(n-1)^2}{4}+\alpha^2,
\end{equation}  
where $\alpha=\alpha(r,l)$ is the smallest positive solution to the equation $\overline S_{l}(\alpha,r)=0.$ 
As in the odd-dimensional case, we may assume that $\alpha r\to \Phi$ as $r\to \infty$ for some $\Phi>0;$ we are going to prove first that 
\begin{equation}\label{alpha-konv-2}
\alpha\sim\frac{\pi}{r}\ \ {\rm as}\ \ r\to \infty.
\end{equation}

By the integral representation of ${\bf P}_{-\frac{1}{2}+i\gamma }$ (see Robin \cite{Robin}) we have that 
\begin{eqnarray*}
	\overline S_{1}(\gamma,x) &=& {\bf P}_{-\frac{1}{2}+i\gamma }(\cosh(x))=\frac{\sqrt{2}}{\pi}\coth(\gamma\pi)\int_x^\infty\frac{\sin(\gamma z)}{\sqrt{\cosh(z)-\cosh(x)}}{\rm d}z\\&=&
	\frac{\sqrt{2}}{\pi}\coth(\gamma\pi)\int_0^\infty\frac{\sin(\gamma(x+y))}{\sqrt{\cosh(x+y)-\cosh(x)}}{\rm d}y\\&=&\frac{\sqrt{2}}{\pi}\coth(\gamma\pi)\left[\cos(\gamma x)\int_0^\infty\frac{\sin(\gamma y)}{\sqrt{\cosh(x+y)-\cosh(x)}}{\rm d}y\right.\\&&\ \ \ \ \ \ \ \ \ \ \ \ \ \ \  \ \ \   \ \ \ \ \ \ \ \left.+\sin(\gamma x)\int_0^\infty\frac{\cos(\gamma y)}{\sqrt{\cosh(x+y)-\cosh(x)}}{\rm d}y\right].
\end{eqnarray*}
Since $\cosh(x+y)\sim \frac{e^{x+y}}{2}$ and $\cosh(x)\sim \frac{e^{x}}{2}$ as $x\to \infty$, it turns out that 
\begin{equation}
\overline S_{1}(\gamma,x)\sim \frac{{2}}{\pi}\coth(\gamma\pi)\left[p_1(\gamma)\cos(\gamma x)+q_1(\gamma)\sin(\gamma x)\right]e^{-\frac{x}{2}}\ \ {\rm as} \ \ x\to \infty,
\end{equation}
where $$p_1(\gamma)=\int_0^\infty\frac{\sin(\gamma y)}{\sqrt{e^y-1}}{\rm d}y\ \ {\rm and}\ \ q_1(\gamma)=\int_0^\infty\frac{\cos(\gamma y)}{\sqrt{e^y-1}}{\rm d}y,$$
see also Robin \cite[p. 55]{Robin}.
Lebesgue's dominated convergence theorem implies that 
\begin{equation}\label{Lebesgue}
\lim_{\gamma \to 0}p_1(\gamma)=0\ \ {\rm and}\ \ \lim_{\gamma \to 0}q_1(\gamma)=\int_0^\infty\frac{1}{\sqrt{e^y-1}}{\rm d}y=\pi.
\end{equation}
By the stability property of differentiation  with respect to the approximation of hyperbolic functions, the solution $\alpha>0$ of   $\overline S_{l}(\alpha,r)=0$  will be approximated by the smallest positive root $\alpha$ of the equation $S_l^\#(\alpha,r)=0$, where 
$$S_1^\#(\gamma,x)=\left[p_1(\gamma)\cos(\gamma x)+q_1(\gamma)\sin(\gamma x)\right]e^{-\frac{x}{2}}\ \ {\rm and}\ \  S_{k}^\#(\gamma,x)={\frac{\partial S_{k-1}^\#}{\partial x}(\gamma,x)}e^{-x}, $$
for every $k\geq 2,$ $\gamma,x>0.$
%
%
%
%
Accordingly, $S_l^\#(\alpha,r)=0$ is equivalent to 
\begin{equation}\label{egyenlet-paros-3}
p_l(\alpha)\cos(\alpha r)+q_l(\alpha)\sin(\alpha r)=0,
\end{equation}
where 
\begin{eqnarray}\label{P-Q-2}
\left\{
\begin{array}{lll}
p_{k+1}(\gamma)&=&\gamma q_k(\gamma)+(\frac{1}{2}-k)p_k(\gamma), \\
q_{k+1}(\gamma)&=&-\gamma p_k(\gamma)+(\frac{1}{2}-k)q_k(\gamma),
\end{array}\right. k\geq 1,\ \gamma> 0.
\end{eqnarray}
In particular, relations (\ref{P-Q-2}) and (\ref{Lebesgue}) imply that 
\begin{equation}\label{Lebesgue02}
\lim_{\gamma \to 0}p_k(\gamma)=0\ \ {\rm and}\ \ \lim_{\gamma \to 0}q_k(\gamma)=(-1)^{k-1}\pi\frac{(2k-3)!!}{2^{k-1}},\ \ k\geq 2.
\end{equation}
Taking the limit $r\to \infty$ (thus $\alpha\to 0$) in  (\ref{egyenlet-paros-3}), the latter limits give  that $\sin(\Phi)=0$, i.e., $\Phi=\pi,$ which concludes the proof of (\ref{alpha-konv-2}). 

We now determine $c_l\in \mathbb R$ such that 
\begin{equation}\label{alpha-ujabb-3}
\alpha\sim\frac{\pi}{r}+\frac{c_l}{r^2}\ \ {\rm as}\ \ r\to \infty.
\end{equation}
By (\ref{egyenlet-paros-3}) one has approximately 
$p_l\left(\frac{\pi}{r}+\frac{c_l}{r^2}\right)\cos\left(\frac{c_l}{r}\right)+q_l\left(\frac{\pi}{r}+\frac{c_l}{r^2}\right)\sin\left(\frac{c_l}{r}\right)=0,$
thus
$$c_l=-\pi\lim_{z\to 0}\frac{p_l(z)}{zq_l(z)}=-\pi \frac{p_l'(0)}{\lim_{z \to 0}q_l(z)}.$$   
The first relation of (\ref{P-Q-2}) implies the recurrence relation $p_{k+1}'(0)=\lim_{z \to 0}q_k(z)-\frac{1}{2}p_k'(0),$  $k\geq 1$, with 
$$p_1'(0)=\int_0^\infty\frac{y}{\sqrt{e^y-1}}{\rm d}y=2\pi\ln 2.$$
Consequently, by (\ref{Lebesgue}) and (\ref{Lebesgue02}), one has $c_1=-2\pi\ln 2$,  and  $$c_l=2\pi\left(1+\frac{1}{3}+...+\frac{1}{2l-3}-\ln 2\right),\ \ l\geq 2,$$
which completes the proof. \hfill $\square$

\begin{remark}\rm \label{remark-n-5} (i) Let $n=5$ $($thus $l=2)$. The transcendental equation $S_2(\alpha,r)=0$ is equivalent to  (\ref{S2-egyenlet}). If $r\to 0$, a similar reasoning as above shows that $\alpha\to \infty$ and $\alpha \sim \frac{\Phi}{r}$ as $r\to 0$ for some $\Phi> 0.$ Thus, taking the limit in $\alpha r\cos(\alpha r) \frac{\tanh(r)}{r}-\sin(\alpha r)=0$ as $r\to 0$, it yields  $\Phi \cos(\Phi)-\sin(\Phi)=0$, which is equivalent to $J_{\frac{3}{2}}(\Phi)=0.$ Since $\Phi> 0$ is minimal with the latter property, it follows that $\Phi=j_{\frac{3}{2},1}$. Now, for some $c_0\in \mathbb R$ let $\alpha^2 \sim \frac{\Phi^2}{r^2}+c_0$ as $r\to 0$. Again by (\ref{S2-egyenlet}) and using the fact that $\tan(\Phi)=\Phi$,  a Taylor expansion implies that $c_0=-\frac{2}{3}.$ Thus, (\ref{elso-sajat-alpha}) provides 
	$$\lambda_{1,d_h}(B_r)=4+\alpha^2\sim 4+\frac{\Phi^2}{r^2}+c_0=\frac{j_{\frac{3}{2},1}^2}{r^2}+\frac{10}{3}\ \ {\rm as}\ r\to 0,$$ which is exactly (\ref{nullaban-aszimptotikus}) for $n=5$. A similar argument applies in  higher odd-dimensions as well. 
	
	(ii)	On the right hand side of relations (\ref{alpha-ujabb})  and (\ref{alpha-ujabb-3}) the exponent $2$ cannot be replaced by any other number $s\in \mathbb R$; if $s<2$ then $c_l=0$, while if $s>2$ then $|c_l|=\infty.$
\end{remark}

\medskip

\section{Proof of Theorem \ref{fotetel-CD}}\label{section-Cheng}

Let $\kappa,r>0$ and the integer $n\geq 2$ be fixed. 
By the  proof of Theorem \ref{fotetel-hiperbolikus} we recall that 
$$
\int_{B_r^\kappa}|\nabla_{g_h}w^*|_{g_h}^2{\rm d}v_{g_h}=\lambda_{1,d_h}(B_r^\kappa)\int_{B_r^\kappa}(w^*)^2{\rm d}v_{g_h},
$$
where $w^*$ is from (\ref{w-csillag}).  
By the latter relation and the differentiation formula (\ref{F-differential}) one has the identity
\begin{equation}\label{hiperbolikus-diff-egyenloseg}
\int_{B_r^\kappa}R_{n+2}^2(d_{h}(0,x))\sinh^2\left({\kappa d_{h}(0,x)}\right){\rm d}v_{g_h}(x)=\frac{\kappa^2 n^2}{\lambda_{1,d_h}(B_r^\kappa)}\int_{B_r^\kappa}R_{n}^2(d_{h}(0,x)){\rm d}v_{g_h}(x),
\end{equation}
the function  $R_\theta:[0,r]\to \mathbb R$ being defined by
 $$R_\theta(\rho)={\bf F}\left(\frac{\theta-1}{2}+ i\alpha^\kappa_{n,r},\frac{\theta-1}{2}-i\alpha^\kappa_{n,r};\frac{\theta}{2};-\sinh^2\left(\frac{\kappa \rho}{2}\right)\right),\ \theta,\rho>0,$$
where $\alpha^\kappa_{n,r}$ is given in (\ref{alpha-ertek}). 

In the sequel, we summarize those properties of $R_\theta$ which will play crucial roles in our proof. 

\begin{proposition}\label{R-tulajdonsagok} The following properties hold: 
	\begin{itemize}
		\item[(i)] $\rho \mapsto R_n(\rho)$ is positive and decreasing on $[0,r)$ with $R_n(0)=1$ and $R_n(r)=0;$
			\item[(ii)] $\rho \mapsto R_{n+2}(\rho)$ is positive on $[0,r];$
			\item[(iii)] $\rho \mapsto \frac{R_n(\rho)}{R_{n+2}(\rho)\sinh(\kappa \rho)}$ is decreasing on $(0,r].$
	\end{itemize} 
\end{proposition}

{\it Proof.} (i) 
We notice that $R_n(\rho)=f(\tanh(\kappa \rho/2))$, where $f$ is from (\ref{meg-ez-is-kell}). Thus  $R_n(0)=1$ and since $\rho=r$ is the first positive solution to the equation $R_n(\rho)=0$, see (\ref{meg-ez-is-kell}) and (\ref{utolso-az-elso-tetelnel}), the function  $\rho \mapsto R_n(\rho)$ is positive  on $[0,r)$. Moreover, by (\ref{ODE}) one can easily see that $f$ is decreasing on $\left[0,\tanh(\frac{\kappa r}{2})\right)$, so is 
$R_n$ on $[0,r)$. 

(ii) $R_{n+2}(0)=1$ and the differentiation formula (\ref{F-differential})  yields $$R_{n+2}(\rho)=-\frac{\kappa n}{\lambda_{1,d_h}(B_r^\kappa)\sinh(\kappa \rho)}R_{n}'(\rho)>0,\ \rho\in (0,r].$$

(iii) By the continued fraction representation (15.7.5) of \cite{Digital}, it turns out that  $$\frac{R_n(\rho)}{R_{n+2}(\rho)\sinh(\kappa \rho)}=T(\coth(\kappa \rho)),$$ where  $$\ds T(t)=\ds x_0t-\frac{y_1}{x_1t-\ds\frac{y_2}{x_2t-\ds\frac{y_3}{\ddots}}},\ t>0,$$
with $x_l=\frac{n+2l}{2}$  and $ y_l=\frac{1}{4}\left(l^2+l(n-1)+\frac{\lambda_{1,d_h}(B_r^\kappa)}{\kappa^2}\right),\ \l\geq 0.$ Since $T$ is increasing and $\coth$ decreasing on $[0,\infty)$, the proof is complete.  \hfill $\square$\\

Due to Propositions \ref{layer-cake} and \ref{R-tulajdonsagok}/(i), it follows that
\begin{equation}\label{jobb-oldal}
\int_{B_r^\kappa}R_{n}^2(d_{h}(0,x)){\rm d}v_{g_h}(x)=\int_{0}^rR_{n}^2(\rho)\sinh^{n-1}\left({\kappa \rho}\right){\rm d}\rho.
\end{equation}
On the other hand, since $\rho\mapsto R_{n+2}^2(\rho)\sinh^2(\kappa \rho)$ is a $BV$-function, we can represent it as the difference of two decreasing functions $q^1$ and $q^2$. Thus, one can apply  Proposition \ref{layer-cake} for the  functions $\rho\mapsto q^i(\rho)-q^i(r)$  $(i=1,2)$, obtaining that 
$$\int_{B_r^\kappa}R_{n+2}^2(d_{h}(0,x))\sinh^2\left({\kappa d_{h}(0,x)}\right){\rm d}v_{g_h}(x)=\int_0^rR_{n+2}^2(\rho)\sinh^{n+1}\left({\kappa \rho}\right){\rm d}\rho.$$
Accordingly, by  (\ref{hiperbolikus-diff-egyenloseg}) and (\ref{jobb-oldal}), we obtain the identity
\begin{equation}\label{ujabb-azonossag}
{\lambda_{1,d_h}(B_r^\kappa)}\int_0^r R_{n+2}^2(\rho)\sinh^{n+1}\left({\kappa \rho}\right){\rm d}\rho={\kappa^2 n^2}\int_{0}^rR_{n}^2(\rho)\sinh^{n-1}\left({\kappa \rho}\right){\rm d}\rho.
\end{equation}

\subsection{Proof of (\ref{2-becsles})} 
Assume the contrary of (\ref{2-becsles}), i.e., 
$\lambda_{1,d}(B_r(x_0))> \lambda_{1,d_h}(B_r^\kappa).$ Taking $\delta_0>0$ sufficiently small, one has 
\begin{equation}\label{Faber-Krahn-anti}
\int_{B_r(x_0)} |\nabla u|_{{{d}}}^2{\rm d}\mu>(\lambda_{1,d_h}(B_r^\kappa)+\delta_0)\int_{B_r(x_0)} u^2{\rm d}\mu, \ \ \forall u\in {\rm Lip}_0(B_r(x_0)).
\end{equation}
 Let 
 \begin{equation}\label{w-definicio}
 w(x)=f\left(\tanh\left(\frac{\kappa d(x_0,x)}{2}\right)\right)\equiv R_n(d(x_0,x)), \ \ x\in B_r(x_0),
 \end{equation}
 where $f$ is from (\ref{meg-ez-is-kell}).  Due to (\ref{utolso-az-elso-tetelnel}), one has  $w(x)=0$ for every $x\in \partial B_r(x_0)$. Moreover, by using elementary truncations, one can construct a sequence of nonnegative functions $w_k\in {\rm Lip}_0(B_r(x_0))$ such that $\{w_k\}_k$ and $\{|\nabla w_k|_d\}_k$  converge pointwisely  to $w$ and $|\nabla w|_d$ in $B_r(x_0)$, respectively, and by the properness of $(M,d)$ (i.e., every bounded
 and closed subset of $M$ is compact), the support of both $w_k$ and $|\nabla w_k|_d$ is the compact set $\overline{B_{r-\frac{1}{k}}(x_0)}$ for every $k\in \mathbb N$.  Applying (\ref{Faber-Krahn-anti}) for $w_k$, the Lebesgue dominated convergence theorem implies that $w$  verifies 
 \begin{equation}\label{Faber-Krahn-anti-2}
 \int_{B_r(x_0)} |\nabla w|_{{{d}}}^2{\rm d}\mu\geq \left(\lambda_{1,d_h}(B_r^\kappa)+\delta_0\right)\int_{\Omega} w^2{\rm d}\mu> \lambda_{1,d_h}(B_r^\kappa)\int_{B_r(x_0)} w^2{\rm d}\mu.
 \end{equation}
 Relation (\ref{w-definicio}), the non-smooth chain rule  and  the eikonal inequality  $|\nabla d(x_0,\cdot)|_d(x)\leq 1$ for every $x\in M\setminus\{x_0\}$ imply  that  
$$|\nabla w|_{{{d}}}(x)=\left|R_n'(d(x_0,x))\right|\cdot|\nabla d(x_0,\cdot)|_d(x)\leq R_{n+2}(d(x_0,x)) \frac{\lambda_{1,d_h}(B_r^\kappa)\sinh(\kappa d(x_0,x))}{\kappa n}.$$
Therefore, a similar reasoning as in (\ref{hiperbolikus-diff-egyenloseg})-(\ref{ujabb-azonossag}),  Proposition \ref{layer-cake} and (\ref{Faber-Krahn-anti-2}) give that 
\begin{equation}\label{ujabb-egyenlotlenseg}
{\lambda_{1,d_h}(B_r^\kappa)}\int_0^r R_{n+2}^2(\rho)\sinh^{2}\left({\kappa \rho}\right){\sf A}^\mu_{\rho}(x_0){\rm d}\rho>{\kappa^2 n^2}\int_{0}^rR_{n}^2(\rho){\sf A}^\mu_{\rho}(x_0){\rm d}\rho.
\end{equation}
 For further use, let $\Psi:(0,r]\to \mathbb R$ be the continuous function defined by 
\begin{equation}\label{g-definicio}
\Psi({\rho})=1-\frac{\kappa^{n-1}}{n\omega_n}\frac{{\sf A}^\mu_{\rho}(x_0)}{{\sinh}^{n-1}(\kappa \rho)}.
\end{equation}
Hypothesis $({\bf BG})^{n,\kappa}_{x_0}$  implies that $\Psi$ is non-decreasing on $(0,r).$ Moreover, by  the local density assumption {\rm $({\bf D})^n_{x_0}$} it turns out that $\limsup_{{\rho}\to 0}{\Psi({\rho})}=0,$
thus $\Psi$ is non-negative; accordingly,  one has  
 $${\sf A}^\mu_{\rho}(x_0)\leq \frac{n\omega_n}{\kappa^{n-1}}{\sinh}^{n-1}(\kappa \rho),\ \rho>0.$$ 
In particular, the latter inequality shows that the terms in (\ref{ujabb-egyenlotlenseg}) are well defined.

For the sake of simplicity, we introduce the function $H:[0,r]\to \mathbb R$ defined by
$$H({\rho}):={\lambda_{1,d_h}(B_r^\kappa)} R_{n+2}^2(\rho)\sinh^{2}\left({\kappa \rho}\right)-{\kappa^2 n^2}R_{n}^2(\rho).$$ With this notation the identity (\ref{ujabb-azonossag}) and inequality (\ref{ujabb-egyenlotlenseg}) can be rewritten into  the forms
\begin{equation}\label{H-iden}
\displaystyle\int_0^rH({\rho})\sinh^{n-1}\left({\kappa \rho}\right){\rm d}{\rho} =0\ \ {\rm and}\ \  \displaystyle\int_0^rH({\rho}){\sf A}^\mu_{\rho}(x_0){\rm d}{\rho} >0,
\end{equation}
respectively. 
Since 
$$\lim_{{\rho}\to 0}\frac{R^2_n(\rho)}{R^2_{n+2}(\rho)\sinh^2(\kappa \rho)}=+\infty\ \ {\rm and}\ \ \lim_{{\rho}\to r}\frac{R^2_n(\rho)}{R^2_{n+2}(\rho)\sinh^2(\kappa \rho)}=0,$$  Proposition \ref{R-tulajdonsagok}/(iii) implies that the equation $H({\rho})=0$ has a unique solution in $(0,r)$; let us denote by ${\rho}_0\in (0,r)$ this element. The above arguments also show that
\begin{equation}\label{rho-tulaj}
H({\rho})<0,\ \forall {\rho}\in [0,{\rho}_0)\ \ {\rm and}\ \  H({\rho})>0,\ \forall {\rho}\in ({\rho}_0,r].
\end{equation}
%
 Since ${\sf A}^\mu_{\rho}(x_0)=\frac{n\omega_n}{\kappa^{n-1}}(1-\Psi({\rho})){\sinh}^{n-1}(\kappa \rho),$ see (\ref{g-definicio}), the two relations of (\ref{H-iden}) imply 
\begin{equation}\label{H-egyenlotlenseg}
\displaystyle\int_0^r \Psi(\rho){\sinh}^{n-1}(\kappa \rho) H({\rho}){\rm d}{\rho} < 0.
\end{equation}
By relations (\ref{H-egyenlotlenseg}), (\ref{rho-tulaj}), (\ref{H-iden})   and the monotonicity of $\Psi$ we have 
\begin{eqnarray*}
	0&>& \displaystyle\int_0^r \Psi(\rho){\sinh}^{n-1}(\kappa \rho) H({\rho}){\rm d}{\rho}=\int_0^{\rho_0} \Psi(\rho){\sinh}^{n-1}(\kappa \rho) H({\rho}){\rm d}{\rho}+\int_{\rho_0}^r \Psi(\rho){\sinh}^{n-1}(\kappa \rho) H({\rho}){\rm d}{\rho}\\&\geq &\Psi(\rho_0)\int_0^{\rho_0} {\sinh}^{n-1}(\kappa \rho) H({\rho}){\rm d}{\rho}+\int_{\rho_0}^r \Psi(\rho){\sinh}^{n-1}(\kappa \rho) H({\rho}){\rm d}{\rho}  \\&= &-\Psi(\rho_0)\int_{\rho_0}^r {\sinh}^{n-1}(\kappa \rho) H({\rho}){\rm d}{\rho}+\int_{\rho_0}^r \Psi(\rho){\sinh}^{n-1}(\kappa \rho) H({\rho}){\rm d}{\rho} \\&= &\int_{\rho_0}^r \left[-\Psi(\rho_0)+\Psi(\rho)\right]{\sinh}^{n-1}(\kappa \rho) H({\rho}){\rm d}{\rho}\\&\geq&0,
\end{eqnarray*}
a contradiction, which implies the validity of the inequality (\ref{2-becsles}). 

\subsection{Equality in (\ref{2-becsles})} Assume we have equality in (\ref{2-becsles}), i.e.,  $$\lambda_{1,d}(B_r(x_0))= \lambda_{1,d_h}(B_r^\kappa).$$ In particular, the latter relation implies that $$\int_{B_r(x_0)} |\nabla w|_{{{d}}}^2{\rm d}\mu\geq  \lambda_{1,d_h}(B_r^\kappa)\int_{B_r(x_0)} w^2{\rm d}\mu,$$
where $w$ is defined in (\ref{w-definicio}). By a similar reasoning as in the previous part we arrive (instead of (\ref{H-egyenlotlenseg})) to the inequality
\begin{equation}\label{H-egyenlotlenseg-2}
\displaystyle\int_0^r \Psi(\rho){\sinh}^{n-1}(\kappa \rho) H({\rho}){\rm d}{\rho} \leq  0.
\end{equation}
Accordingly, a similar estimate as above -- based on (\ref{H-egyenlotlenseg-2}) --  implies that
$$0\geq \int_{\rho_0}^r \left[-\Psi(\rho_0)+\Psi(\rho)\right]{\sinh}^{n-1}(\kappa \rho) H({\rho}){\rm d}{\rho}\geq 0,$$
thus $\Psi(\rho)=\Psi(\rho_0)$ for every ${\rho}\in (\rho_0,r).$ Having this relation,  we similarly have by (\ref{rho-tulaj}) and the monotonicity of $\Psi$ that 
\begin{eqnarray*}
	0&\geq& \displaystyle\int_0^r \Psi(\rho){\sinh}^{n-1}(\kappa \rho) H({\rho}){\rm d}{\rho}=\int_0^{\rho_0} \Psi(\rho){\sinh}^{n-1}(\kappa \rho) H({\rho}){\rm d}{\rho}+\int_{\rho_0}^r \Psi(\rho){\sinh}^{n-1}(\kappa \rho) H({\rho}){\rm d}{\rho} \\&= &\int_0^{\rho_0} \Psi(\rho){\sinh}^{n-1}(\kappa \rho) H({\rho}){\rm d}{\rho}+\Psi(\rho_0)\int_{\rho_0}^r {\sinh}^{n-1}(\kappa \rho) H({\rho}){\rm d}{\rho}  \\&= &\int_0^{\rho_0} \left[\Psi(\rho)-\Psi(\rho_0)\right]{\sinh}^{n-1}(\kappa \rho) H({\rho}){\rm d}{\rho}\\&\geq &0.
\end{eqnarray*}
Thus, we have $\Psi(\rho)=\Psi(\rho_0)$ for every $\rho\in (0,\rho_0).$ Summing up, we have $\Psi(\rho)=\Psi(\rho_0)$ for every $\rho\in (0,r).$ Since $\limsup_{{\rho}\to 0}{\Psi({\rho})}=0$, it follows that $\Psi\equiv 0$ on $(0,r)$.  By (\ref{g-definicio}) it turns out that
$$\mu(B_\rho(x_0))=\int_0^\rho {\sf A}^\mu_{t}(x_0){\rm d}t=n\omega_n\int_0^\rho \left(\frac{\sinh(\kappa t)}{\kappa}\right)^{n-1}{\rm d}t=V_\rho^\kappa,\ \ \ \rho\in (0,r),$$
which concludes the proof.  \hfill $\square$

\medskip

\section{Proof of Theorem \ref{fotetel-Finsler}}\label{section-Funk}
For simplicity, let $M:=B^n=\{x\in \mathbb R^n:|x|<1\}$ be the $n$-dimensional Euclidean unit ball, $n\geq 2$, and consider the \textit{Funk metric} $F:B^n\times
\mathbb R^{n}\to \mathbb R$  defined by
\begin{equation}\label{F-a-metrika}
F(x,y)=\frac{\sqrt{|y|^2-(|x|^2|y|^2-\langle
		x,y\rangle^2)}}{1-|x|^2}+\frac{\langle x,y\rangle}{1-|x|^2},\ x\in
B^n,\ y\in T_xB^n=\mathbb R^n.
\end{equation}
Hereafter, $|\cdot|$ and
$\langle\cdot, \cdot\rangle$ denote the $n$-dimensional Euclidean
norm and inner product. The pair $(B^n,F)$ is a non-reversible Finsler manifold which falls into the class of Randers spaces,  see Cheng and Shen
\cite{Cheng-Shen} and Shen \cite{Shen-konyv}. The co-metric of $F$ is
\begin{equation}\label{F-a-polar-0}
F^*(x,y)=|y|-\langle
	x,y\rangle
,\ \ \ (x,y)\in B^n\times
\mathbb R^{n}.
\end{equation}
The distance function associated
to $F$ is 
$$d_{F}(x_1,x_2)=\ln\frac{\sqrt{|x_1-x_2|^2-(|x_1|^2|x_2|^2-\langle x_1,x_2\rangle^2)}-\langle x_1,x_2-x_1\rangle}{\sqrt{|x_1-x_2|^2-(|x_1|^2|x_2|^2-\langle x_1,x_2\rangle^2)}-\langle x_2,x_2-x_1\rangle},\ x_1,x_2\in B^n,$$
see Shen \cite[p.141 and p.4]{Shen-konyv}; in particular,
$$d_{F}(0,x)=-\ln(1-|x|) \ {\rm and}\ \ d_{F}(x,0)=\ln(1+|x|),\ \ x\in B^n,$$ thus $(B^n,d_F)$ is a quasimetric space. 
The Busemann-Hausdorff volume form on $(B^n,F)$ is ${\text d}v_F(x)={\rm d}x,$ see Shen \cite[Example 2.2.4]{Shen-konyv}. The Finsler manifold $(B^n,F)$ is forward (but not backward) complete, it has constant negative flag curvature ${\bf K}= -\frac{1}{4}$, see Shen \cite[Example 9.2.1]{Shen-konyv}, and constant negative weighted Ricci curvature $\mathrm{Ric}_N(v)=-\frac{n-1}{4}
-\frac{(n+1)^2}{4({N-n)}}$ for $N\in (n,\infty),$ $\mathrm{Ric}_\infty(v)=-\frac{n-1}{4}$ and $\mathrm{Ric}_n(v)=-\infty$ for every $v\in T_xB^n$ with $F(x,v)=1,$  see Ohta \cite{Ohta-Pacific}. 
In particular,  $(B^n,F)$ does not satisfy the ${\sf CD}(-(n-1)\kappa^2,n)$ condition for any $\kappa>0.$  

\subsection{First proof of Theorem \ref{fotetel-Finsler} (via Theorem \ref{fotetel-CD})} 
We observe that  
$$\mu(B_\rho(0))=\int_{B_\rho(0)}{\text d}v_F(x)=\int_{|x|<1-e^{-\rho}}{\text d}x=\omega_n(1-e^{-\rho})^n,\ \ \rho>0,$$
therefore,  
$${\sf A}^\mu_{\rho}(0)=n\omega_n (1-e^{-\rho})^{n-1}e^{-\rho},\ \ \rho>0.$$ Thus,  
it is easy to prove that {\rm $({\bf D})^n_{0}$} holds and the function $\rho\mapsto \frac{{\sf A}^\mu_{\rho}(0)}{\sinh^{n-1}(\kappa \rho)}$ $(\rho>0)$ is decreasing for every $\kappa>0$; thus $({\bf BG})^{n,\kappa}_{0}$ holds for \textit{every} $\kappa>0$. By applying Theorem \ref{fotetel-CD}, it follows that	$$\lambda_{1,d_F}(B_r(0))\leq \lambda_{1,d_h}(B_r^\kappa)$$
for \textit{every}  $r,\kappa>0$. 
Taking $r\to \infty$ and $\kappa\to 0$ in the latter inequality,  Savo's estimate (\ref{SAvo}) (or Theorem \ref{fotetel-hiperbolikus}) yields that  $$\lambda_{1,d_F}(B^n)= \lim_{r\to \infty} \lambda_{1,d_F}(B_r(0))\leq 0,$$ which concludes the proof.\hfill $\square$


\subsection{Second proof of Theorem \ref{fotetel-Finsler} (direct estimate of the fundamental frequency)} \label{sub-5-2}
By (\ref{osszefugges-F-gradiens}) and definition  (\ref{first-eigenvalue-general}), it turns out that 
\begin{equation}\label{elso-Finsler}
\lambda_{1,d_F}(B^n)=\inf_{u\in H_{0,F}^1(B^n)\setminus \{0\}}\frac{\ds\int_{B^n} F^{*2}(x,Du(x)){\rm d}v_F(x)}{\ds\int_{B^n} u^2(x){\rm d}v_F(x)},
\end{equation}
where  $H_{0,F}^{1}(B^n)$ is the closure of $C_0^\infty(B^n)$
with respect to the (positively homogeneous) norm
\begin{equation}\label{Sobolev-norm}
\|u\|_{F}=\left(\displaystyle \int_{B^n} F^{*2}(x,Du(x)){\rm
	d}v_{F}(x)+\displaystyle \int_{B^n} u^2(x){\rm
	d}v_{F}(x)\right)^{1/2},
\end{equation}
see Ge and Shen \cite{GS}, Ohta and Sturm \cite{OS}. 

For every $\alpha>0$, let 
\begin{equation}\label{profil-fct}
u_\alpha(x):=-e^{-\alpha d_F(0,x)}=-(1-|x|)^\alpha,\ \ x\in B^n.
\end{equation}
By (\ref{tavolsag-derivalt}) or by a direct computation we have  $F^*(x,Du_\alpha(x))=\alpha(1-|x|)^\alpha,$ thus
$$\ds\int_{B^n} F^{*2}(x,Du_\alpha(x)){\rm d}v_F(x)=\alpha^2\ds\int_{B^n} (1-|x|)^{2\alpha}{\rm d}x=\alpha^2n\omega_n {\sf B}(2\alpha+1,n),$$
where ${\sf B}$ denotes the Beta function. In a similar way, one has 
$$\ds\int_{B^n} u^2_\alpha(x){\rm d}v_F(x)=n\omega_n {\sf B}(2\alpha+1,n).$$
Accordingly, $u_\alpha\in H_{0,F}^1(B^n)$ for \textit{every} $\alpha>0$; thus, the functions $u_\alpha$ can be used as test functions in (\ref{elso-Finsler}), obtaining that  
$$\lambda_{1,d_F}(B^n)\leq \inf_{\alpha>0}\frac{\ds\int_{B^n} F^{*2}(x,Du_\alpha(x)){\rm d}v_F(x)}{\ds\int_{B^n} u^2_\alpha(x){\rm d}v_F(x)}=\inf_{\alpha>0}\alpha^2=0,$$
which ends the proof.

\subsection{Third proof of Theorem \ref{fotetel-Finsler} (via the Finsler-Laplace operator)} \label{sub-5-3} Finally, we provide a moral explanation of the fact that $\lambda_{1,d_F}(B^n)=0$  by using directly the Finsler-Laplace operator. More precisely, for every $0<\rho<1,$ we consider the eigenvalue problem 
\begin{equation}\label{laplace-finsler-egyenlet}
\left \{
\begin{array}{lll}
-\Delta_{F}w=\lambda_{\rho} w & {\rm in} &  B_\rho^n, \\
w=0 & {\rm on} &  \partial B_\rho^n,%
\end{array}
\right.
\end{equation}
where  $B_\rho^n=\{x\in \mathbb R^n:|x|<\rho\}$ and $\lambda_{\rho}>0$. 
Having in our mind the shape of the function in (\ref{profil-fct}), we look for the eigenfunction in (\ref{laplace-finsler-egyenlet}) in the form $w(x)=f(|x|)$ for some $f:[0,\rho]\to \mathbb R$ enough smooth, verifying also $f\leq 0$ and $f'\geq 0$ on $(0,\rho)$. 
If such a function $w:=w_\rho$ is not zero, we clearly have that 
\begin{equation}\label{minek-utolso}
\lambda_{1,d_F}(B_\rho^n)\leq \lambda_{\rho}.
\end{equation}


One has $$Dw(x)=f'(|x|)\frac{x}{|x|},\ x\neq 0.$$  Moreover, due to (\ref{F-a-polar-0}), the Legendre transform is given by $$J^*(x,y):=\frac{\partial}{\partial y}\left(\frac{1}{2}F^{*2}(x,y)\right)= (|y|-\langle
x,y\rangle)\left(\frac{y}{|y|}-x\right);$$ thus,  
$${\nabla}_F w(x)=J^*(x,Dw(x))=f'(|x|)(1-|x|)^2\frac{x}{|x|},\ x\neq 0.$$
Accordingly, since $\sigma_F(x)=1$ (due to ${\text d}v_F(x)={\rm d}x$),   it turns out that
\begin{eqnarray*}
	\Delta_Fw(x)&=&{\rm div}_F({\nabla}_F w(x))={\rm div}({\nabla}_F w(x))\\&=&\frac{\rm d}{{\rm d}s}(f'(s)(1-s)^2)\big|_{s=|x|}+f'(|x|)(1-|x|)^2\frac{n-1}{|x|},\ x\neq 0.
\end{eqnarray*}
The latter computations together with (\ref{laplace-finsler-egyenlet}) give the second order ODE
\begin{equation}\label{ODE-utolso}
(f'(s)(1-s)^2)'+f'(s)(1-s)^2\frac{n-1}{s}+\lambda_{\rho} f(s)=0,\ s\in (0,\rho),
\end{equation}
subject to the boundary  condition $f(\rho)=0.$ Standard ODE arguments show that the non-singular solution at the origin for (\ref{ODE-utolso}) is  given by means of the Gaussian hypergeometric function as
\begin{eqnarray*}
	f(s)&=&-(1-s)^{-\frac{1+\sqrt{1-4\lambda_\rho}}{2}}\times\\&&\times {\bf F}\left(\frac{n-1}{2}+\frac{\sqrt{n^2-4\lambda_\rho}-\sqrt{1-4\lambda_\rho}}{2},\frac{n-1}{2}-\frac{\sqrt{n^2-4\lambda_\rho}+\sqrt{1-4\lambda_\rho}}{2};n-1;s\right),
\end{eqnarray*}
for $s\in (0,\rho)$ whenever $\lambda_\rho<\frac{1}{4}$.  The condition $f(\rho)=0$ implies that 
$\lambda_\rho>0$ is the smallest solution of 
$${\bf F}\left(\frac{n-1}{2}+\frac{\sqrt{n^2-4\lambda_\rho}-\sqrt{1-4\lambda_\rho}}{2},\frac{n-1}{2}-\frac{\sqrt{n^2-4\lambda_\rho}+\sqrt{1-4\lambda_\rho}}{2};n-1;\rho\right)=0.$$
Since ${\bf F}\left(n-1,-1;n-1;1\right)=0$, a continuity argument shows that $\lambda_\rho\to 0$ as $\rho\to 1.$ Therefore, by (\ref{minek-utolso})  we have	
$$\lambda_{1,d_F}(B^n)= \lim_{\rho\to 1}\lambda_{1,d_F}(B_\rho^n) \leq \lim_{\rho\to 1}\lambda_{\rho}=0,$$
which concludes the proof.
 \hfill $\square$

\begin{remark}\rm (i)
	 The symmetrization of the Funk distance $d_F$ is the  \textit{Klein distance} on $B^n$,   given by  $$d_K(x_1,x_2)=\frac{1}{2}\left(d_F(x_1,x_2)+d_F(x_2,x_1)\right),\ \ x_1,x_2\in B^n,$$
	while its corresponding  Klein metric on $B^n$ is 
	$$F_K(x,y)=\frac{\sqrt{|y|^2-(|x|^2|y|^2-\langle
			x,y\rangle^2)}}{1-|x|^2},\ x\in
	B^n,\ y\in T_xB^n=\mathbb R^n,$$ see Cheng and Shen
	\cite{Cheng-Shen} and Shen \cite{Shen-konyv}.  The Klein volume form is  ${\text
		d}v_{K}(x)={(1-|x|^2)^{-\frac{n+1}{2}}}{\text d}x.$ Clearly, $(B^n,F_K)$ is a complete Riemannian manifold with constant negative  sectional curvature $-1$. According to (\ref{elso-hatarertek}), we have 
	\begin{equation}\label{Klein-concord}
	\lambda_{1,d_K}(B^n)=\frac{(n-1)^2}{4}.
	\end{equation}
If we mimic the argument from \S\ref{sub-5-2} on the Klein model $(B^n,F_K)$ for functions of the form $$w_\gamma(x)=e^{-\gamma d_K(0,x)},\ \ x\in B^n,$$ we should assume  that $\gamma>\frac{n-1}{2}.$ In fact, $w_\gamma$ belongs to $H_{0,F_K}^1(B^n)$ if and only if $\gamma>\frac{n-1}{2};$ 
indeed, the integrals $$\ds\int_{B^n} w^2_\gamma(x){\rm d}v_K(x)= \int_{B^n}\left(\frac{1-|x|}{1+|x|}\right)^\gamma {(1-|x|^2)^{-\frac{n+1}{2}}}{\text d}x=n\omega_n\int_0^1(1-t)^{\gamma-\frac{n+1}{2}}(1+t)^{-\gamma-\frac{n+1}{2}}t^{n-1}{\rm d}t $$
and 
$$\int_{B^n} F_K^{*2}(x,Dw_\gamma(x)){\rm d}v_K(x)=\gamma^2\int_{B^n} w^2_\gamma(x){\rm d}v_K(x)$$
 are convergent if and only if $\gamma-\frac{n+1}{2}>1$. The above  computations show that 
$$ \lambda_{1,d_K}(B^n)\leq \inf_{\gamma>\frac{n-1}{2}}\frac{\ds\int_{B^n} F_K^{*2}(x,Dw_\gamma(x)){\rm d}v_K(x)}{\ds\int_{B^n} w^2_\gamma(x){\rm d}v_K(x)}=\inf_{\gamma>\frac{n-1}{2}}\gamma^2=\frac{(n-1)^2}{4},$$
which together with McKean's lower estimate (\ref{Mckeannn})  provides a perfect concordance with (\ref{Klein-concord}). 

(ii) We notice that  $H_{0,F_K}^1(B^n)$ is the usual Sobolev space over the Riemannian Klein model $(B^n,F_K)$, see e.g. Hebey  \cite{Hebey}. However, the set $H_{0,F}^1(B^n)$ over the non-reversible Finslerian Funk model  $(B^n,F)$ endowed with the norm (\ref{Sobolev-norm}) is not even a vector space, see Krist\'aly and Rudas \cite{Kristaly-Rudas}. Indeed, although  $u_\alpha\in H_{0,F}^1(B^n)$ for every $\alpha>0$, it turns out that  $-u_\alpha\notin H_{0,F}^1(B^n)$ for any $\alpha\in (0,\frac{1}{2}]$,  since
$$\ds\int_{B^n} F^{*2}(x,-Du_\alpha(x)){\rm d}v_F(x)=\alpha^2\int_{B^n} (1+|x|)^2({1-|x|})^{2\alpha-2}{\rm d}x=+\infty.$$ 
\end{remark}

\vspace{0.5cm}
\noindent {\bf Acknowledgment.} The author is grateful to  Denis Borisov and Pedro Freitas for the conversations concerning their paper \cite{BF}. He  also thanks  \'Arp\'ad Baricz, Csaba Farkas, Mihai Mih\u ailescu  and Tibor Pog\'any for their help in special functions and eigenvalue problems. \\

\end{document}